\documentstyle[bezier,11pt]{article}
 \title{  }

\author{  }
\date{Mars 2000}

\catcode`\Ÿ=\active
\defŸ{\"u}

\catcode`\Ž=\active
\defŽ{\'e}

\catcode`\‰=\active
\def‰{\^a}
\catcode`\ˆ=\active
\defˆ{\`a}
\catcode`\=\active
\def{\c c}
\catcode`\=\active
\def{\^e}
\catcode`\=\active
\def{\`e}
\catcode`\•=\active
\def•{\"\i }
\catcode`\"=\active
\def"{\^\i }
\catcode`\™=\active
\def™{\^o}
\catcode`\ž=\active
\defž{\^u}
\catcode`\=\active
\def{\`u}
\catcode`\¶=\active
\def¶{\delta }
\catcode`\¹=\active
\def¹{\pi }

\newcommand{\ds}{\displaystyle}
\newcommand{\fdem}{$\Box$}

\renewcommand{\L}{{I\!\!L}}   
   
\newcommand{\HH}{{I\!\!H}}       
\newcommand{\N}{{I\!\!N}}     
\newcommand{\Z}{{Z\!\!\!Z}}
\newcommand{\R}{{I\!\!R}}     
\newcommand{\D}{{I\!\!D}}
\newcommand{\SS}{\mbox{\rm S\hspace{-0.8ex}\rule[0.05ex]{0.1ex}{1.5ex}
\hspace{0.1ex}}}    
\newcommand{\C}{\mbox{\rm C\hspace{-1.1ex}\rule[0.3ex]{0.1ex}{1.2ex}
\hspace{1.1ex}}}

\newcommand{\cC}{{\cal C}}
\newcommand{\cH}{{\cal H}}
\newcommand{\cB}{{\cal B}}

\newcommand{\cD}{{\cal D}}

\newcommand{\tC}{{\tilde C}}

\def\md{{\buildrel{{\scriptscriptstyle\Delta}}\over \times }}

\begin{document}
\centerline{\bf \Large  Mesures de Hausdorff}
\vspace{5mm}
\centerline{\bf \Large  
de l'ensemble limite de groupes 
kleiniens }
\vspace{5mm}
\centerline{\bf \Large  
g\'eom\'etriquement finis}

\vspace{1cm}

\noindent {\it Nous nous proposons ici de montrer que l'exposant de 
Poincar\'e d'un groupe g\'eom\'etriquement fini co{\"\i}ncide avec la 
dimension de Hausdorff de son ensemble limite et de comparer les 
mesures naturelles port\'ees par cet ensemble : la mesure de 
Patterson et les mesures de Hausdorff et packing pour la jauge 
standart.}

\vspace{1cm}

\noindent {\bf I. Introduction et notations}

\vspace{5mm}

Soit $\HH^{3} = \C\times ]0, +\infty[$ le demi-espace sup\'erieur
de dimension $3$
muni de la m\'etrique hyperbolique 
$\ds{ds^{2}= \frac{{dx^{2}+dy^{2}}}{w^{2}}}$ au 
point $(x+iy, w)$.
On consid\`erera aussi le mod\`ele de la boule $\D^{3}
= \{x \in \R^{3}/ \vert \vert x \vert \vert <1\}$
o\`u cette fois on pose $\ds{ds^{2}= \frac{{\vert dx^{2}\vert 
}}{(1-\vert \vert x\vert \vert^{2})^{2}}}$. Dans les deux cas on note 
$d$ la distance hyperbolique sur $\HH^{3}$. Le bord \`a l'infini 
$\partial 
\HH^{3}$  s'identifie, selon le mod\`ele choisi,  \`a
$\C$ ou \`a $\SS^{2}$ ; dans ce dernier cas, le bord   
$\partial 
\HH^{3}$ est muni de la m\'etrique euclidienne $\vert. \vert$.
Dans la suite on fixe une origine $o$ dans $\HH^{3}$ ; c'est le 
centre de la boule $\D^{3}$ dans le mod\`ele du disque et le 
point $(0, 1)$ dans le mod\`ele du demi-espace.

\noindent Le groupe  des isom\'etries de $\HH^3$ est 
le groupe de Mobius de 
dimension $2$. Une isom\'etrie
$g$ de $\HH^{3}$
agit sur le bord 
 $\partial 
\HH^{3}$
de l'espace hyperbolique par transformation conforme,
le coefficient de conformit\'e en un point $\xi$ \'etant 
$\vert g'(\xi)\vert =
\exp (-\cB_{\xi}(\gamma^{-1}.o, o))$ o\`u
$\cB_{\xi}(., .)$ d\'esigne le cocycle de Buseman au point $\xi$ d\'efini 
par
$$
\forall x, y \in \HH^{3} 
\quad 
\cB_{\xi}(x, y)
= \lim_{z \to \xi} d(x, z)-d(y, z).
$$
De plus, pour tout $\xi, \eta \in \partial \HH^{3}$ on a
$$\vert g.\xi-g.\eta\vert^{2} = \vert g'(\xi)\vert\vert g'(\eta)\vert 
\vert \xi-\eta\vert^{2}.$$
On consid\`ere dans ce qui suit  un groupe kleinien, 
c'est-\`a-dire une sous-groupe discret $G$ de $PSL(2, \C)$. Un tel groupe 
agit proprement discontinuement sur $\HH^{3}$.
{\it L'ensemble limite } $\Lambda_{G}$ de $G$ est 
\'egal \`a $\overline{G.x}- G.x$ ;
c'est le plus petit ferm\'e $G$-invariant de $\partial 
\HH^{3}$ et il  ne d\'epend pas du 
point $x$. L'ensemble des  points $\xi \in \Lambda_{G}$ pour lesquels
il existe un voisinage born\'e du rayon g\'eod\'esique 
$[o, \xi)$ contenant une infinit\'e de points de $G.o$ joue un r\^ole 
crucial dans ce qui suit ; cet ensemble est appel\'e {\it ensemble 
limite radial}, il est $G$-invariant 
et on le note 
$\Lambda_{G}^r$ .
Si l'on pose
$G = \{g_{n}, n \geq 1\}$ on peut \'ecrire 
$$
\Lambda_{G}^r= \bigcup_{k\in \N}
\limsup_{n \to +\infty}O_{x}(g_{n}.x, k)
$$
o\`u  pour tous points
$x, y$ de $\HH^{3}$ et tout $k>0$ on a 
$O_{x}(y, k)=\{\eta \in \partial \HH^{3}/d([x, \eta), y) <k\}$ (on dit 
que $O_{x}(y, k)$ est {\it l'ombre} sur 
$\partial \HH^{3}$
vue de $x$  de la boule de centre $y$ et de rayon 
$k$).
  
\noindent Le groupe $G$ est dit {\it \'el\'ementaire} lorsque
$\Lambda_{G}$ est r\'eduit \`a un ou deux  points ; dans les autres 
cas le cardinal de $\Lambda_{G}$ est infini.   Lorsque $\sharp \Lambda_{G}= 
2$, le groupe $G$ est dit {\it hyperbolique} , il est cyclique 
et engendr\'e par une isom\'etrie hyperbolique. Lorsque
$\sharp \Lambda_{G}= 1$, le groupe $G$ est {\it 
parabolique} et isomorphe soit \`a $\Z$ soit \`a $\Z^2$ 
selon que,  dans le mod\`ele du demi-espace, le groupe $G$ est engendr\'e par 
une ou deux translations ind\'ependantes.

\noindent On  associe \`a $G$ la{\it s\'erie  de 
 Poincar\'e} ${\cal P}_{G}$ d\'efinie , pour tout r\'eel $s\geq 0$ et tous 
points $x, y \in \HH^{3}$,  par
$$
{\cal P}_{G}(x, y, s) = \sum_{g \in G}e^{-sd(x, g.y)}.
$$
L'estimation du volume des boules de $\HH^{3}$ permet de montrer que 
${\cal P}_{G}(x, y, s)<+\infty$ pour tout $s>2$.
On note $\delta_{G}$ (ou encore $\delta$ lorsqu'il n'y a pas 
d'ambiguit\'e) l'exposant critique de cette s\'erie :
$$
\delta _{G}= \inf\{s\geq 0/{\cal P}_{G}(x, y, s) <+\infty\}.$$ Cet exposant
ne d\'epend pas du choix de 
$x$ et $y$ et est appel\'e {\it exposant de Poincar\'e } de $G$.
Le groupe $G$ est dit {\it divergent} lorsque $ {\cal P}_{G}(x, y, 
\delta_{G}) = +\infty$ et {\it convergent} sinon.
Nous avons  le r\'esultat \'el\'ementaire suivant :

\noindent {\bf Lemme 1.1 -}
{\it 
Pour tout $\alpha \in \C^{*} $ l'exposant de Poincar\'e du groupe $P_{\alpha}$
engendr\'e
par la transformation parabolique  $z\mapsto z+\alpha$ est \'egal \`a 
$1/2$.

\noindent Pour tout $\alpha, \beta \in \C^{*} $  ind\'ependant sur $\R$,
l'exposant de Poincar\'e du groupe 
$P_{\alpha, \beta}$ engendr\'e
par les transformations paraboliques  $z\mapsto z+\alpha$ et $z\mapsto 
z+\beta $ est \'egal \`a 
$1$.

\noindent Dans les deux cas les groupes paraboliques sont divergents. 
De plus, si  $G $ est un groupe non \'el\'ementaire contenant un 
 sous-groupe parabolique de rang $k$ alors $\delta>k/2$.
}

\noindent {\it D\'emonstration-}  Dans les deux cas, le point fixe du 
groupe parabolique consid\'er\'e est 
$+\infty$ dans le mod\`ele du demi-espace. Rappelons que si $(z, w)$ 
et $(z', w')$ sont deux points de $\HH^{3}$ leur distance est donn\'ee
par
$$d((z, w), (z', w')) = \log \frac{1+t}{1-t} \quad \mbox{\rm avec}
\quad 
t = \sqrt{\frac{\vert z-z'\vert^{2}+\vert y-y'\vert^{2}}{
\vert z-z'\vert^{2}+\vert y+y'\vert^{2}}}. $$

\noindent Si $p_{\alpha} $ d\'esigne la tranlation $z \mapsto z+\alpha$ on a 
$$d((i, 1), p_{\alpha}^n(i, 1))= d((i, 1), (i+n\alpha, 1)) =
(\log n^{2}) (1+\epsilon(n))$$ et la s\'erie de Poincar\'e de 
$P_{\alpha}$ se comporte comme la s\'erie $\ds{\sum_{n \in \Z^{*}}
\frac{1}{\vert n \vert^{2s}}}$ ; l'exposant de Poincar\'e de $P_{\alpha}$ est 
donc $1/2$ et le groupe est de type
divergent.

\noindent De m\^eme, pour tous entiers $n, m \in \Z^{*}$ on a 
$d((i, 1), p_{\alpha}^np_{\beta}^m(i, 1))= 
(\log \vert n \alpha + m \beta \vert ^{2}) (1+\epsilon(n, m))$
 et la s\'erie de Poincar\'e de 
$P_{\alpha, \beta}$ se comporte comme la s\'erie  double 
$\ds{\sum_{n, m  \in \Z^{*}}\frac{1}{\vert n \alpha +m \beta \vert ^{2s}}}$ ; 
ainsi 
l'exposant de Poincar\'e de $P_{\alpha, \beta}$ est 
 $1$ et ce groupe est  divergent.

\noindent Enfin, si $G$ est un groupe  non \'el\'ementaire  contenant un 
sous-groupe parabolique $P$ on a de fa\c{c}on 
\'evidente $\delta\geq \delta_{P}$. L'in\'egalit\'e stricte est 
plus d\'elicate \`a obtenir et 
d\'ecoule du fait que $G$ est   non \'el\'ementaire  et que les 
groupes paraboliques sont de type divergent. En  effet
si  $G$ est non \'el\'ementaire, il contient  une isom\'etrie 
hyperbolique
 $h$ dont les points  fixes  sont distincts de 
 $+\infty$.  Quitte \`a remplacer $h$ par une puissance suffisamment  
 grande et $P$ par un sous-groupe parabolique de m\^eme rang, 
 on peut  supposer que $G$ contient le produit libre 
 $P*<h>$ ; ainsi 
 $\{h^{n_{1}}p_{1}h^{n_{2}}p_{i_{2}}\cdots
 h^n_{l}p_{l}/ n_{j}\in \Z^{*}, p_{j} \in P-\{Id\}, l\geq 1\}$ est contenu 
 dans $G$ et il vient
 \begin{eqnarray*}
     \sum_{g\in G}e^{-sd(x, g.x)}&\geq &
     \sum_{l\geq 1}\sum_{p_{1}\cdots p_{l} }\sum_{n_{1}, \cdots, 
     n_{l} }e^{-sd(x, h^{n_{1}}p_{1}\cdots
 h^{n_{l}}p_{l}.x)}
     \\
 &\geq & \sum_{l\geq 1}\sum_{p_{1}\cdots p_{l} }\sum_{n_{1}, \cdots, 
     n_{l} }
     e^{-sd(x, h^{n_{1}}.x)}e^{-sd(x, p_{1}.x)}\cdots
 e^{-sd(x, h^{n_{l}}x)}e^{-sd(x, p_{l}.x)}\\
 &\geq & \sum_{l\geq 1}\Bigl( \sum_{p_\in P-\{Id\}}e^{-sd(x, 
 p.x)}
 \sum_{n \in \Z^{*}}
      e^{-sd(x, h^{n}.x)}\Bigr)^l.\\
 \end{eqnarray*}
Le groupe $P$ \'etant divergent, on  a
$\ds{\lim_{\stackrel{s\to \delta_{P}}{s>\delta_{P}}}
\sum_{p_\in P-\{Id\}}e^{-sd(x, 
 p.x)}=+\infty}$ et l'on peut donc choisir  $s_{0} > \delta_{P}$ tel que
$$ 
\sum_{p_\in P-\{Id\}}e^{-s_{0}d(x, 
 p.x)}
 \sum_{n \in \Z^{*}}
     \Bigl( e^{-s_{0}d(x, h^{n}.x)} > 1.$$
Ainsi ${\cal P}_{G}(x, x, s_{0})$ diverge d'o\`u $\delta_{G}> \delta _{P}$.

\noindent Dans ce qui suit le groupe $G$ sera suppos\'e non 
\'el\'ementaire.
On note $C (\Lambda_{G})$ l'enveloppe convexe de 
$\Lambda_{G}$; cet ensemble est invariant sous  l'action de 
$G$ et le quotient $N(G)= C (\Lambda_{G})/G$ est 
le coeur de Nielsen de la vari\'et\'e $M(G) = \HH^{3}/G$. 
On notera $N_{\epsilon}(G)$ un $\epsilon$-voisinage
 de $N(G)$. Lorsque $\Lambda(G) = \Lambda_{G}^r
$, l'ensemble   $N(G)$
est relativement compact ; on dit que $G$ est {\it convexe co-compact}. 
L'\'etude des groupes  non  co-compacts mais de 
 co-volume fini  (c'est-\`a-dire tels 
que la vari\'et\'e $M(G)$ est non compacte mais de volume fini) est 
une premi\`ere \'etape dans l'etude de groupe plus g\'en\'eraux. 
Ces groupes font partie de la classe plus large des groupes dits
{\it g\'eom\'etriquement finis} : 

\noindent {\bf D\'efinition 1.2 -}
 {\it On dit qu'un groupe kleinien non \'el\'ementaire $G$ est g\'eom\'etriquement 
 fini si et seulement si
 il existe $\epsilon >0$ tel que
 $N_{\epsilon}$ est de volume fini.}

\noindent La perte de compacit\'e de $N_{\epsilon}$ se traduit par 
l'appparition dans $\Lambda_{G}$ de points limites non coniques : les 
points   {\it paraboliques born\'ees}.
 
 \noindent {\bf D\'efinition 1.3 -}
 {\it Soit $G$ un groupe kleinien non \'el\'ementaire. Un point $\xi 
 \in \Lambda_{G}$ est dit parabolique born\'e si son stabilisateur 
 $P$ dans $G$ est un groupe parabolique et si 
 $\Lambda_{G}-\{\xi\}$ admet un domaine fondamental relativement 
 compact pour l'action de $P$}.
 
 \noindent Rappelons que  la 
 finitude g\'eom\'etrique  peut \^etre caract\'eris\'ee de 
fa\c{c}on \'equivalent comme suit :
  
 -  pour tout $\epsilon >0$ le volume de $N_{\epsilon}(G)$ est fini.
 
 - pour tout $\epsilon >0$, la partie $\epsilon$-\'epaisse 
 $N(G)^{>\epsilon}$ est relativement compacte.
 
 - le groupe $G$ contient un nombre fini de classes de conjugaison  de 
 groupes paraboliques dont les points fixes  
 sont born\'es.
 
 - le coeur de Nielsen $N(G)$ peut se d\'ecomposer en $C_{0}\cup 
 C_{1} \cdots \cup C_{l}$ o\`u $C_0$ est un ensemble relativement 
 compact et o\`u, pour chaque $i = 1, \cdots, l,$
 il existe un groupe  parabolique $P_{i} \subset G$
et une 
 horoboule  ${\cal H}_{i}$ bas\'ee en $\xi_i$ tels que 
 $C_{i}$ soit isom\'etrique au 
 quotient de 
 ${\cal H}_i\cap C(\Lambda(G))$ par le groupe $P_{i}$
 (notons  que le point fixe $\xi_{i}$  de $P_{i}$ est alors 
 n\'ecessairement born\'e, que le groupe $P_{i}$ agit sur $C(\Lambda_{G})\cap
 \partial \cH_{i}$ o\`u $\partial \cH_{i}$ d\'esigne l'horisph\`ere 
 qui
 borde l'horiboule  $ \cH_{i}$ et que cette action admet  un domaine 
 fondamental relativement compact).

 \noindent {\bf 2. Mesure conforme sur l'ensemble limite d'un groupe 
 kleinien}
 
 \noindent Nous nous int\'eressons \`a la structure de $\Lambda_{G}$ 
 d'un point de vue de la th\'eorie de la mesure. Rappelons tout 
 d'abord la d\'efinition suivante
 
 \noindent {\bf D\'efinition 2.1 -}
 {\it Une mesure finie $\sigma$ sur $\SS^{2}$ est dite $G$-conforme 
 d'exposant $\alpha \in \R$ si, pour tout $g \in  G$,  on a }
 $$
 \frac{d(g^{*}\sigma)}{d\sigma)}(\xi) = \exp (\alpha 
 \cB_{\xi}(\gamma^{-1}.o, o)).
 $$
 
 \noindent Rappelons le proc\'ed\'e de Patterson permettant de 
 construire une mesure $G$-conforme sur $\SS^{2}$ d'exposant 
 $\delta_{G}$. Pour chaque $s>\delta_{G}$ et chaque point $x \in \HH^{3}$
 on note $\sigma_{x}^s$ la mesure orbitale
 $$
 \sigma_{x}^s= \frac{1}{{\cal P}_{G}(o, s)}
\sum_{g \in G}\exp (-s d(x, g.o)) D_{g.o}$$
o\`u $D_{g.o}$ d\'esigne la masse de Dirac en $g.o$.
Lorsque le groupe $G$ est de  type divergent,
toute valeur d'adh\'erence (pour la topologie de la convergence 
\'etroite)
de la famille 
$(\sigma_{x}^s)_{x, s}$ est port\'ee par
$\Lambda_{G}$; on peut alors  montrer que lorsque $s \to \delta_{G}$ 
par valeurs sup\'erieures
 la famille de mesures $(\sigma_{x}^s)_{s}$  converge \'etroitement 
 vers une  mesure $\sigma_{x}$
port\'ee par $\Lambda_{G}$ et v\'erifiant les deux conditions suivantes
$$
\sigma_{x'}(.) = \exp (-\delta\cB_{.}(x', x)) \sigma_{x}(.)
\quad \mbox{\rm et} \quad g^{*}\sigma_{x} = \sigma_{g^{-1}.x} 
$$
o\`u $g^{*}\sigma_{x} $ est la mesure sur $\SS^d$ d\'efinie par 
$g^{*}\sigma_{x}(B) =  \sigma_{x}(g B)$ pour tout bor\'elien $B $ 
de $\SS^d$. 
On dit que la famille $(\sigma_{x})_{x \in \HH^{3}}$ est une 
{\it densit\'e $G$-conforme d'exposant $\delta_{G}$}.

\noindent Il est difficile a priori de montrer qu'un groupe $G$ est 
de type divergent; en particulier D. Sullivan a \'etabli cette 
propri\'et\'e pour les groupes g\'eom\'etriquement finis, en 
\'etudiant le type des densit\'es $\delta_{G}$-conformes de ces groupes. 
Pour ce faire, il faut pouvoir construire de telles densit\'es  
$\delta_{G}$-conformes, , sans savoir \`a priori si $G$ est de type 
convergent ou 
divergent  ;  en
utilisant un argument du \`a Patterson, on modifie
l\'eg\`erement la s\'erie de Poincar\'e en posant
$${\cal P}'_{G}(x, y, s) = \sum_{g \in G}e^{-sd(x, g.y)}h(d(x, g.y))$$
o\`u $h$ est une fonction croissante de $\R^{+}$ dans $\R^{+}$ telle que
les s\'eries ${\cal P}_{G}(x, y, s) $ et ${\cal P}'_{G}(x, y, s)$ aient le m\^eme 
exposant critique et
$$
\forall \eta >0,  \exists  t_{\eta}>0,  \forall t \geq t_{\eta},  \forall s 
\geq 0 \quad
 h(t+s) \leq h(t) e^{\eta s}.
$$

\noindent  Int\'eressons nous maintenant aux propri\'et\'es locales 
 d'une densit\'e $G$-conforme. 
Nous avons le 

\noindent {\bf  Lemme 2.2} (-Th\'eor\`eme de l'ombre de Sullivan-) {\it
Soit $G$ un groupe non \'el\'ementaire et $\sigma $ une densit\'e 
$G$-conforme d'exposant $\alpha$. Il existe $C>1$ et $r_{0}>0$ tel 
que pour tout $r\geq r_{0}$ et tout $g \in G$ on ait}
$$
\frac{1}{C} e^{-\alpha d(o, g.o)} \leq  
\sigma_{x}(O_{x}(g.o, r))\leq C  e^{-\alpha d(o, g.o)+2\alpha }.
$$

\noindent Nous donnerons dans le paragraphe 5 une d\'emonstration d'une 
version un peu plus pr\'ecise de ce lemme. Soulignons  en cependant 
deux cons\'equences  importantes    :

- en remarquant qu'une ombre
$O_{x}(g.o, r)$ rencontre au plus un nombre uniform\'ement born\'e 
d'ombres 
$O_{x}(h.o, r)$ avec $h \in G$ et $d(o, g.o) -1\leq d(o, h.o) \leq 
d(o, g.o) +1$, on montre gr\^ace \`a ce lemme que  l'existence d'une
densit\'e $G$-conforme $\sigma$ d'exposant $\alpha$ 
entra\^ine 
$$
\sharp\{ g \in G/ d(o, g.o) \leq R\} \leq C e^{\alpha R}
$$
et donc  $\alpha \geq \delta_{G}$ par d\'efinition de 
l'exposant  de Poincar\'e de $G$. Il n'existe donc pas de densit\'e 
$G$-conforme d'exposant  $< \delta_{G}$.

- si $\ds{\sum_{n\geq 1}e^{-\alpha d(o, g.o)}<+\infty }$
(ce qui est le cas lorsque $\alpha >\delta$) alors toute densit\'e
$\alpha$-conforme donne une mesure nulle \`a $\Lambda_{G}^r$.

\noindent {\bf Th\'eor\`eme 2.3 -}
{\it Si $G$ est g\'eom\'etriquement fini alors pour tout $x \in 
\HH^{3}$ on a
$\sigma_{x}(\Lambda_{G}^c)= \sigma_{x}(\Lambda_{G})$ et le groupe $G$ est de type 
divergent.}

\noindent {\it D\'emonstration-}
Il suffit de d\'emontrer 
que si $\xi $ est un point parabolique born\'e de $\Lambda_{G}$ alors 
il est de $\sigma$-mesure nulle.
Rappelons que pour tout voisinage ouvert $V$ de $\xi$ dans 
$\HH^{3}\cup \partial \HH^{3}$ on a
$$
\sigma_{x}(\xi)\leq 
\sigma_{x}(V) \leq \liminf_{s_{i} \to \delta^{+}}
\sigma_{x}^{s_i}(V).
$$
Il suffit donc d'exhiber des ouverts $V$  contenant $\xi$ et tels que
$\ds{\liminf_{s_{i} \to \delta^{+}}
\sigma_{x}^{s_i}(V)}$ soit arbitrairement petit.
Nous allons nous placer dans le mod\`ele du demi-espace et poser 
$\xi = +\infty$. On consid\`ere un domaine fondamental $\cD_{\infty}$
pour l'action 
de $P = stab_{G}(\xi)$ sur $\partial \HH^{3} - \{+\infty\}$
 que l'on peut choisir de fa\c{c}on que $\Lambda_{G} \cap \cD_{\infty} 
 $ 
 soit 
relativement compact dans
$\partial \HH^{3} - \{+\infty\}$. On note $\cD$ le c\^one  sur 
$\cD_{\infty}$
issu de 
$+\infty$ et $G'$ le sous-ensemble de $G$ contenant les isom\'etries 
$g$ telles que $g.o \in \cD$; les valeurs d'adh\'erence de $G'.o$ sont 
donc 
contenues dans $\overline{\cD_{\infty}}$.

\noindent Choisissons $x \in \HH^{3}$ tel que $B_{\infty}(x, 
g'.o)<0$ pour tout $g \in G'$. Notons $P = \{p_k, k \geq 1\}$ et 
choisissons  une suite d\'ecroissante d'ouverts
$(V_{n})_{n \geq 1}$ tels que  pour tout $n$ on ait
$\ds{
G.o\cap V_{n} \subset \cup_{k\geq n}p_{k}\cD .
}$ Il vient imm\'ediatement
$$
\sigma_{x}^{s_{i}}(V_{n})\leq \frac{1}{{\cal P'}(o, s_{i})}
\sum_{k\geq n}\sum_{g' \in G'}e^{-s_{i}d(x, p_{k} g'.o)}
$$
La convexit\'e des horisph\`eres  et le choix de $x$ font que l'angle 
au point $x$ entre les segments $[x, p_{k}.x]$ et $[x, g'.o]$ est 
minor\'e par une constante strictement positive ; il existe donc une 
(autre) constante  $C>0$ telle 
que 
$d(x, p_{k}g'.o) \geq d(x, p_{k}x)+ d(x, g'.o)-K$
gr\^ace au fait suivant

\noindent {\bf Lemme 2.4 -}
{\it Pour tout $\epsilon>0$ il existe $C_{\epsilon}>0$
 tel que pour tout triangle g\'eod\'esique  de $\HH^{3}$
 de c\^ot\'es $a, b, c$ et d'angles oppos\'es $\alpha, \beta$ et $\gamma$ 
 avec $\gamma \geq \epsilon$ on a} 
 $$a+b-C_{\epsilon}\leq c\leq a+b.$$
{\it D\'emonstration-} La loi du cosinus en g\'eom\'etrie hyperbolique
donne $$cosh\ c\  = \ cosh\  a \ cosh\  b\  - \ sinh\  a \ sinh\  b\  cos
\ \gamma 
\geq cosh\  a\  cosh\  b \ (1-\vert cos \ \gamma \vert )$$
si bien que
$$c \geq \log\  cosh \ c \geq \log (cosh\  a\  cosh\  b(1-\vert cos\  \gamma 
\vert ) )\geq a+b -C_{\gamma}$$
avec $C_{\gamma} = 2 \log 2 -\log (1-\vert cos  \gamma \vert ) \to +\infty$ 
lorsque $\gamma \to 0$.\fdem 

\noindent Prenons alors $\eta >0$ tel que $\delta_{G}> \delta_{P} +\eta$ et 
choisissons $x \in \HH^{3}$ pour que 
$d(x, g.o)\geq t_{\eta}$ ; il vient
$$\sigma_{x}^{s_{i}}(V_{n}) \leq e^{K s_{i}}\sum_{k \geq n}
e^{-(s_{i}-\eta) d(x, p_{k}.x)} \frac{1}{{\cal P}_{G}'(s_{i})}\sum_{g' \in G'}
h(x, g'.o)e^{-s_{i}d(x, g'.o)}$$
d'o\`u l'on d\'eduit imm\'ediatement
$$
\liminf_{s_{i}\to \delta}\sigma_{x}^{s_{i}}(V_{n})
\leq 
e^{K \delta}\sum_{k \geq n}
e^{-(\delta-\eta) d(x, p_{k}.x)} \sigma_{x}(\Lambda_{G}).
$$
Comme $\delta-\eta > \delta_{P}$ on a 
$\ds{\lim_{n \to +\infty} \sum_{k \geq n} e^{-(\delta-\eta)d(x, 
p_{k}.x)}=0}$ d'o\`u le r\'esultat escompt\'e.
Ainsi $\sigma_{x}$ ne charge pas $\xi$. L'ensemble des points paraboliques \'etant 
au plus d\'enombrable on a  $\sigma(\Lambda_{G}^r) = 
\sigma(\Lambda_{G})$ et un argument de type Borel-Cantelli permet  
de conclure que $G$ est de type divergent. \fdem

\noindent {\bf 3 - Mesure de Bowen-Margulis des groupes g\'eom\'etriquement 
finis} 

\noindent Rappelons le proc\'ed\'e de  Sullivan qui permet d'associer \`a
une mesure
de Patterson $\sigma_x$  une mesure
 sur le fibr\'e unitaire tangent   de la vari\'et\'e $\HH^{3}/G$,  invariante par 
 le flot g\'eod\'esique.

\noindent Notons
$\SS^d \md \SS^d$ l'ensemble $\SS^d\times \SS^d$ priv\'e
de sa diagonale.  On peut identifier le
 fibr\'e unitaire tangent $T^1\HH^{3}$ au produit
$\SS^d \md \SS^d \times \R$
en associant \`a un \'el\'ement $v=(y, \vec{v})\in T^1 /HH^{3}$
le triplet
$(\xi^-, \xi^+, r)$    o\`u
$\xi^-$ et $\xi^+$ sont les extr\'emit\'es de
la g\'eod\'esique orient\'ee d\'etermin\'ee par $(y, \vec{v})$
et $r = B_{\xi^+}(y, o)$. Dans ces coordonn\'ees,
 l'action  d'une isom\'etrie
$g$ de $G$  est donn\'ee par
$$
g(\xi^-, \xi^+, r)= (g \xi^-, g \xi^+,
r+B_{\xi^+}(x,g^{-1} x))
$$
tandis que le flot g\'eod\'esique $(\tilde{\phi}_t)_{t \in \R}$ agit sur
$T^1\HH^{3}$ par
$$
\tilde{\phi}_t(\xi^-, \xi^+, r)= (\xi^-, \xi^+, r-t).
$$
Puisque $\sigma _x$ est $\delta_{G}$-conforme, la mesure
$\ds{\frac{\sigma _x(d\xi^-)\sigma _x(d\xi ^+)}{
\vert \xi^--\xi^{+}\vert ^{2\delta }}}$
est une $G$-invariante sur   $\SS^d \md \SS^d$ :
c'est le {\it courant g\'eod\'esique} $c^{\sigma}$ associ\'e \`a 
$\sigma = (\sigma _x)$.
La mesure $\tilde \mu ^{\sigma}= c^{\sigma}\otimes dt$
est invariante sous les actions de $G$ et du flot g\'eod\'esique
 $(\tilde \phi_t)$;
 son support est   $\Lambda _{G}\md \Lambda _{G} \times \R$.
 Elle induit donc par
passage au quotient
une mesure $\mu ^{\sigma}$, invariante sous l'action du flot g\'eod\'esique
$(\phi_t)$ sur $T^1(M)$ et dont le support est   
$(\Lambda _{G}\md \Lambda _{G}\times \R)/G$.

\noindent Lorsque $G$ est co-compact ou convexe co-compact, 
l'ensemble 
$(\Lambda _{G}\md \Lambda _{G}\times \R)/G$ est compact;
$\mu ^{\sigma}$ est
alors finie et  c'est la mesure d'entropie maximale.
Lorsque $G$ est g\'eom\'etriquement fini, divergent, et contient des
transformations
 paraboliques, la question de la finitude de $\mu
^{\sigma}$ se pose de fa\c{c}on naturelle
(remarquons que si $G$ \'etait convergent, $\sigma _x$ chargerait uniquement
les points
paraboliques et la mesure $\mu ^{\sigma}$ serait clairement de masse infinie).

\noindent {\bf Th\'eor\`eme 3.1.}
{\it Si $G$ est un 
groupe g\'eom\'etriquement fini, alors $\mu^{\sigma}$ est finie.
}

\noindent {\it D\'emonstration -}
La projection sur $\HH^{3}$ du support de $\tilde \mu ^{\sigma}$
 est contenue dans $C(\Lambda _{G})$.
Le groupe $G$ \'etant g\'eom\'etriquement fini, le quotient
$C(\Lambda _{G})/G$ se  d\'ecompose en la r\'eunion
disjointe d'un compact $C_0$
 et d'une famille finie $C_1, \cdots , C_l$  de ``bouts
cuspidaux'': pour $i\geq 1$, chaque $C_i$ est isom\'etrique
au quotient de l'intersection de
 $C(\Lambda _{G})$ et
d'une horiboule $\cH_{\xi _i}$ par un groupe parabolique $P_i$.
Choisissons
  un domaine fondamental bor\'elien $\cC_i$ pour l'action de $G$ sur
la pr\'eimage  de $C_i$ dans
$\HH^{3}$. Sans perte de g\'en\'eralit\'e, on peut supposer que $\cC_0$ est
relativement compact
et que, pour $i\geq 1$,  $\cC_i$ est  un domaine fondamental
pour l'action
de $ P_i$ sur $\cH_{\xi_i}\cap C(\Lambda _{G})$.
On a
 $$
\mu^{\sigma}(T^1(\HH^{3}/G))=\sum_{i=0}^l\tilde \mu^{\sigma}(T^1\cC_i)
 =\sum_{i=0}^l \int_{\SS^d \md \SS^d}
c^{\sigma}(d\xi^- d\xi^+) \int_{(\xi^- \xi^+)\cap \cC_i }dt.
$$
Puisque $\cC_0$ est relativement compact dans $\HH^{3}$, il existe 
$\epsilon>0$ tel que $\vert \xi ^--\xi ^+\vert \geq \epsilon_{0}$
pour toute
g\'eod\'esique $(\xi^-|\xi^+)$ rencontrant
 $\cC_0$. Par cons\'equent,
$\tilde \mu^{\sigma}(T^1\cC_0) \leq \sigma_x(\SS^d)^2
/\epsilon^{2\delta}.$

\noindent Ainsi, la mesure $\mu ^{\sigma}$ est finie si et seulement si
 $\tilde \mu ^{\sigma}(T^1\cC_i)$ est fini
 pour $i=1,\cdots ,l$.
Notons pour simplifier $\cC$ l'un des domaines fondamentaux $\cC_i$,
$P$ le groupe
parabolique   et $\xi$ le point parabolique correspondants.
Puisque $G$ est g\'eo\-m\'e\-tri\-quement fini, on peut choisir
 un domaine fondamental bor\'elien $\cD_{\infty}$
 pour l'action de $ P$ sur $\SS^d -\{\xi\}$
tel que
$  \cD_{\infty}\cap \Lambda_{G} $ soit relativement compact dans
$\SS^d -\{\xi\}$. Le groupe $G$ \'etant divergent, on a
$\sigma_x\{\xi\} = 0$ et donc
$$
\tilde \mu^{\sigma}(T^1\cC)=
\sum_{p, q \in P}\int_{p \cD_{\infty}\times q \cD_{\infty}}
c^{\sigma}(d\xi^- d\xi^+) \int_{(\xi^- \xi^+)\cap \cC}dt.
$$
En utilisant le fait que $c^{\sigma}$ est invariante sous l'action
de $G$, on obtient:
$$
\tilde \mu^{\sigma}(T^1\cC)=
\sum_{p, q \in P}\int_{ \cD_{\infty}\times p^{-1}q \cD_{\infty}}
c^{\sigma}(d\eta^- d\eta^+) \int_{(\eta^- \eta^+)\cap p^{-1}\cC}dt.
$$
Puisque $\cC$ est un domaine fondamental pour l'action de $ P$
sur $\cH_{\xi}\cap C(\Lambda _{G})$, on a donc:
$$
\mu^{\sigma}(T^1\cC)=
\sum_{p \in  P}\int_{ \cD_{\infty}\times p\cD_{\infty}}
c^{\sigma}(d\eta^- d\eta^+) \int_{(\eta^- \eta^+)\cap \cH_{\xi}}dt.
$$
D'un point de vue g\'eom\'etrique, toute g\'eod\'esique
$(\eta^- \eta^+)$ qui passe par $\cH_{\xi}$ se projette sur $M$ en une
g\'eod\'esique
qui fait une incursion dans la r\'egion cuspidale $C$ et le terme
$\int_{(\eta^- \eta^+)\cap \cH_{\xi}}dt$ correspond
 \`a la longueur de cette incursion.
Comme $\cD_{\infty}\cap \Lambda _{G}$
 est relativement compact dans $\SS^d-\{\xi\}$,
 il existe un compact $K \subset \HH^{3}$ contenant $o$ tel que
 toute g\'eod\'esique $(\eta^- \eta^{+})$ issue de 
 $\cD_{\infty}\cap \Lambda _{G}$ et rencontrant $\cH_{\xi}$ 
 traverse $K$ ; une telle g\'eod\'esique v\'erifie donc 
 $\vert \eta^--\eta^{+}\vert \geq \exp( - diam(K))$.
  Si de plus $\eta^+\in p
\cD_{\infty}\cap \Lambda _{G}$,
la g\'eod\'esique $(\eta^+ \eta^-)$ passe par $p(K)$ et la diff\'erence
$|\int_{(\eta^- \eta^+)\cap \cH_{\xi}}dt-d(o, p.o)|$ est donc
major\'ee par $2  diam K$. Finalement
il existe une constante $C>0$
telle que
$$\frac1{C}\sum_{p \in   P} \sigma_o(p \cD_{\infty})(d(o, p.o)-C)
\leq \mu^{\sigma}(T^1\cC)\leq
C\sum_{p \in   P} \sigma_x(p\cD_{\infty})(d(o, p.o)+C).$$
Puisque $\sigma _o$ est $\delta$-conforme, on a  $\sigma_o(p \cD_{\infty})=
\int_{\cD_{\infty}}
e^{-\delta  B_{\eta}(p^{-1} o, o)}\sigma_o(d \eta)$.
Comme
$\cD_{\infty}
\cap \Lambda _{G}$ est relativement compact dans $\SS^d -\{\xi\}$ il existe
$\epsilon>0$
tel que pour tout $\eta  \in \cD_{\infty} \cap \Lambda _{G}$ et
pour tous les $p\in  P$ sauf peut-\^etre un nombre
fini,
l'angle au point  $x$ entre les segments g\'eod\'esiques $[o, p^{-1}o]$ et
$[o, \eta)$ soit sup\'erieur \`a $ \epsilon$ ;  le terme
$\sigma_x(p \cD_{\infty})$ est donc \'equivalent \`a
 $e^{-\delta d(o, po)}$, uniform\'ement en $p\in \  P$. 
 L'in\'egalit\'e 
 $\delta> \delta_{P}$ permet de conclure
 que la s\'erie $\ds{\sum_{p \in   P} d(o, p.o) e^{-\delta d(o, 
 po)}}$ est convergente.
\fdem

\noindent {\bf 4- Dimension de Hausdorff de l'ensemble limite des groupes 
g\'eom\'etriquement finis}

\noindent Nous reprenons ici la d\'emonstration de Sullivan du fait que la 
dimension de Hausdorff de l'ensemble limite 
d'un groupe g\'eom\'etriquement fini est 
\'egale \`a l'exposant critique de ce groupe. Ce r\'esultat est en 
fait valable pour tous les groupes kleiniens non \'el\'ementaires, 
la preuve de Sullivan repose sur la finitude de la mesure de 
Bowen-Margulis.

\noindent Rappelons la d\'efinition de  la dimension de 
Hausdorff
$HD(E)$ 
d'un sous-ensemble 
$E$ d'un espace m\'etrique $X$.
Pour $x\in X $ et  $r>0$ on note $B(x, r)$ la boule ouverte de centre 
$x$ et de rayon $r$. Pour $s$ et $\epsilon$ strictement positifs 
on pose
$$
\mu_{\epsilon}^s(E) = \inf 
\{\sum_{i=1}^nr_{i}^s/ n\geq 1, r_{i}<\epsilon \quad \mbox{et} \quad   E \subset 
\cup_{i=1}^nB(x_{i},r_{i)}\}
$$
puis 
$$\mu^s(E)= \sup_{\epsilon>0}\mu_{\epsilon}^s(E).$$ La dimension de 
Hausdorff de $E$ est alors d\'efinie par
$
\ds{HD(E) = \inf \{s>0/ \mu^s(E) =0 \}.}
$

\noindent L'in\'egalit\'e $HD(\Lambda_{G})\leq \delta$ est 
facile \`a \'etablir. Puisque $\Lambda_{G}$ diff\`ere de 
$\Lambda_{G}^r$ d'au plus une famille d\'enombrable de points, il 
suffit de montrer que $HD(\Lambda_{G}^r) \leq \delta$.
Rappelons que, en \'ecrivant  $G = \{g_{n}/ n \geq 1\}$,  on a
$$
 \Lambda_{G}^r = \bigcup_{k\geq 1} \limsup_{n \to +\infty}O_{o}(g_{n}.o, 
 k).
$$
Fixons $k \geq 1$ et posons  $\ds{\Lambda_k=
\limsup_{n \to +\infty}O_{o}(g_{n}.o, k)} $. L'ombre 
$O_{o}(g.o, k)$ est une boule de $(\SS^d, \vert. \vert)$ de rayon
$\leq C(k) \exp (-d(o, g.o))$ ; pour $\epsilon >0$ fix\'e
on a alors
$$\Lambda_{k} \subset \bigcup_{n/ C(k) \exp (-d(o, 
g.o))<\epsilon}O_{o}(g.o, k).$$ 
Le lemme de l'ombre de Sullivan donne alors 
$\ds{
\mu_{\epsilon}^s(\Lambda_{k})\leq C(k)^s P_{G}(o, o, s)
}.$
Ainsi
$\ds{\mu^{s}(\Lambda_{k})< +\infty}$
si $s> \delta$ ;  cette propri\'et\'e \'etant v\'erifi\'ee pour 
tous les $s> \delta$ on a alors 
$\mu^{s}(\Lambda_{k})=0$ pour tout $s > \delta$ et donc, 
 $\mu_{s}$ \'etant  une mesure, 
$\mu_{s}(\Lambda_{G}^r)=0$. Il vient $HD(\Lambda_{G}) \leq \delta$
par d\'efinition  de la 
dimension de Hausdorff. 

\noindent  L'argument principal pour \'etablir
l'in\'egalit\'e $HD(\Lambda_{G}\geq \delta$ 
repose 
sur un id\'ee de Frostman : si $E$ est un  sous-ensemble de 
$\Lambda_{G}^r$ de mesure de Patterson  strictement positive et dont 
tout recouvrement fini 
par des boules $B_{1}, \cdots, B_{n}$ de rayons respectifs 
$r_{1},  \cdots, r_{n}$ est tel que $\sigma(B_{i}) \leq 
C r_{i}^{\delta(1-\eta)}$ o\`u $\eta$ est une constante 
strictement positive alors
$HD(E) \geq \delta(1-\eta)$. En effet on a alors
$$\sum_{i=1}^nr_{i}^{\delta(1-\eta)}
\geq \frac{1}{C}\sum_{i=1}^n\sigma(B_{i}) \geq \frac{1}{C} \sigma(E)$$
si bien que $E$ a une $\delta(1 -\eta)$-mesure de Hausdorff 
strictement positive et $HD(E) \geq \delta(1 - \eta)$. Il suffit
 pour conclure
 d'exhiber pour chaque valeur de $\eta>0$ un tel ensemble $E$. Pour 
ce faire nous utiliserons les deux propri\'et\'es suivantes :  

\noindent {\bf FLemme 4.1.}
{\it Soit $\pi$ la projection canonique de $T^{1}M$ sur $M$ qui \`a un 
point $v = (x, \vec{v})$ de $T^{1}M$ associe son point de base $\pi(v)=x.$
Si $\mu^\sigma$
 est fini alors }
 $$
 \lim_{t \to +\infty} \frac{1}{t} d(\pi(\phi_{t}v), o) = 0 \quad 
 \mu^{\sigma}(dv)-p.s.
 $$
Posons $\rho(t, v) = d(\pi(\phi_{t}v), o)$ et notons $\rho'$ la 
d\'eriv\'ee de $\rho$ dans la direction du flot $(\phi_t)$ d\'efinie par
$\ds{\rho'(v)= \lim_{s \to 0} \frac{\rho(s, v)-\rho(0, v)}{s}}$. On a 
$\vert \rho'(v)\vert \leq 1$ et donc
$\rho'\in \L^{1}(\mu^\sigma)$ puisque $\mu^\sigma$ est finie. 
Par ailleurs
$\rho'(-v) = -\rho'(v)$ si bien que $\mu^\sigma(\rho') = 0$. Le 
Fait 1 est alors une cons\'equence du th\'eor\`eme 
ergodique de Birkhoff. \fdem

\noindent {\bf Lemme 4.2.} {\it Soit $G$ un groupe de type divergent et
$(\sigma_{x})_{x}$ la densit\'e $\delta$-conforme associ\'ee.
Pour tout $R >0$
il existe $C>0$ tel que pour tout point $y \in \HH^{3}$ on ait}
$$
\sigma_{o}(O_{o}(y, R))\leq C \exp (-\delta d(o, y)+\delta d(y, 
G.o)).
$$
Le groupe $G$ \'etant divergent,   la mesure 
$\sigma_{o}$ est la limite \'etroite lorsque $s \to \delta^{+}$
de la famille de mesures orbitales
$\ds{\sigma_{o}^s = \frac{1}{P_{G}(o, o, s)}
\sum_{g \in G}e^{-s d(o, g.o)}}$.  Notons $\tilde{O}(y, R)$  l'ensemble des 
points $z$ de $\HH^{3}$ tels que le segment $[o, z]$
rencontre la boule de centre $y$ et de rayon $R$ ; si $z \in \tilde{O}(y, R)$,
l'angle entre les 
segments $[o, y]$ et $[y, z]$ est minor\'e par une constante 
strictement positive et d'apr\`es le lemme 2.4. 
il existe  alors   $K>0$ telle
que 
$d(o, z) \geq d(o, y)+d(y, z)-K$ pour tout point $z \in \tilde{O}(y, R)$.
Notons $g_{0}$ l'\'el\'ement de $G$ tel que
$d(y, g_{0}.o) = d(y, G.o)$ ; on a $d(g_{0}.o, z) \leq d(g_{0}.o, 
y)+d(y, z)$. Finalement pour
tout point $z \in \tilde{O}(y, R)$ on a 
$$
d(o, z) \geq d(o, y)+d(g_{0}.o, z)-d(g_{0}.o, y)-K
$$
si bien que 
$$
\sigma_{o}^s(O_{o}(y, R) )\leq
\frac{e^{s(K+d(y, g_{0}.o)-d(o, y))}}{P_{G}(o, o, s)}
\sum_{g/ g.o \in \tilde{O}(y, R)}e^{-sd(g_{0}.o, g.o)}
\leq e^{s(K+d(y, G.o)-d(o, y))} \sigma_{o}^s(\overline{\HH^{3}}).
$$
Il suffit pour conclure de faire tendre $s $ vers 
$\delta.$\fdem

\noindent 
Indiquons maintenant comment D. Sullivan a \'etabli l'in\'egalit\'e 
$HD(\Lambda_{G}) \geq \delta.$ Par le lemme 4.1. et le th\'eor\`eme 
d'Egorov, il existe un compact $V \subset T^{1}M$ de 
$\mu^\sigma$-mesure strictement positive et tel que 
$$\forall \eta >0, \exists t_{\eta}>0, \forall v \in V, \forall t \geq 
t_{\eta} \quad  d(\pi(\phi_{t}.v, o) < \eta t.$$
Notons ${\cal V}$ un domaine fondamental relativement compact pour 
l'action de $G$ sur la pr\'eimage  de $V$ ; sans perdre en 
g\'en'eralit\'e, on peut supposer $o /in {\cal V}$. 
Posons 
$E = \{v(+\infty)/ v \in \tilde{V}\}$ ; on a $\sigma(E)>0$.
Pour tout  $\xi \in E$ notons $\xi_{t}$ le point situ\'e sur
le rayon $[o, \xi)$ \`a distance $t$ de l'origine ; on a
$d(\xi_{t}, G.o) \leq \eta t +diam( \cal{V})$. En utilisant alors le 
lemme 4.2.  il vient
$$
\sigma_{o}(O_{o}(\xi_{t}, 1))\leq C' e^{-\delta t+ \eta \delta t }.
$$
Rappelons alors que $O_{o}(\xi_{t}, 1)$ est une boule de $\SS^d$ 
centr\'ee en $\xi$ et dont le rayon est de l'ordre de $e^{-t}$ ; on 
montre ainsi que toute boule  de rayon 
$r>0$ centr\'ee en un point de $E$  est de mesure de Patterson 
inf\'erieure \`a $C r^{\delta(1-\eta)}$. L'argument de Frostman 
 permet de conclure.

\noindent 
{\bf 5. Comportement local de la mesure de Patterson}

\noindent Dans ce paragraphe nous \'evaluons  la mesure de Patterson des 
boules dont le centre appartient \`a l'ensemble limite  de $G$. Cette 
estimation est due \`a D. Sullivan ; la d\'emonstration qu'il 
proposait s'appuyait sur la notion de limite g\'eom\'etrique des groupes 
Kleinien
mais pr\'esentait une faille quant \`a l'argument final (en effet, il 
est possible de 
construire une suite de groupes g\'eom\'etriquement finis 
convergeant  vers un groupe parabolique, contrairement \`a ce 
qu'annon\c{c}ait D. Sullivan). La d\'emonstration que nous proposons 
ici reprend tout en le simplifiant  l'argument de  B. Stratmann et 
M. Velani ; elle pr\'esente de  plus assez de souplesse pour \^etre 
\' etendue   au cas de la courbure 
variable.

\noindent Dans ce qui suit,  $\xi$ appartient \`a $\Lambda_{G}$  et $\xi(t)$ 
d\'esigne le  point situ\'e sur le segment g\'eod\'esique $[o, \xi)$  
\`a distance $t$ de l'origine. 
On note $V(o, \xi, t)$ l'intersection de $\SS^d$ avec le demi-espace contenant 
$\xi$  orthogonal en $\xi(t)$ au segment  $[o, \xi)$. 
Rappelons que la distance d'un point $x \in \HH^{3}$ \`a une 
g\'eod\'esique $(\xi^-\xi^{+})$ est donn\'ee par
$$
ch\ d(x, (\xi^-\xi^{+}))= 2 \frac{\vert x-\xi^-\vert .\vert x-\xi^+\vert}
{\vert \xi^--\xi^+\vert.(1-\vert x\vert ^{2})}.
$$
On en d\'eduit imm\'ediatement que  
l'ensemble 
$V(o, \xi, t)$ est une boule euclidienne sur $\SS^d$ de centre $\xi$ 
et de rayon $r_{t}$ de l'ordre de $e^{-t}$.

\noindent
Le groupe $G$ \'etant g\'eom\'etriquement fini, le quotient
$C(\Lambda _{G})/G$ se  d\'ecompose en la r\'eunion
disjointe d'un compact $C_0$
 et d'une famille finie $C_1, \cdots , C_l$  de ``bouts
cuspidaux'', l'ensemble  $C_0$ correspondant \`a la partie 
$\epsilon$-\'epaisse du coeur de Nielsen $N(G) $ pour un certain 
$\epsilon >0$. 
La pr\'eimage  de $C_i$ dans
$\HH^{3}$  est not\'ee $\tC_{i}$ et, pour  $1\leq i \leq n,$  on 
d\'esigne par $k_{i}$ le 
rang du point parabolique $\xi_{i}$ correspondant au bout cuspidal 
$C_{i}$.

\noindent Nous nous proposons d'\'etablir le th\'eor\`eme suivant

\noindent {\bf Th\'eor\`eme 5.1.}
{\it Pour tout point $\xi \in \Lambda_{G}$ et 
tout r\'eel $t>0$ on a
$$ 
\sigma_{o}(V(o, \xi, t))\asymp  
e^{ -\delta t+(k_{t}-\delta)d(\xi(t), G.o)}
$$
o\`u $k_{t}= \delta$ si}
$\xi(t) \in \tC_{0}$ et  $k_{t}=k_{i}$  si $\xi(t) \in \tC_{i}.$

\noindent (la notation $a \asymp b $ signifie qu'il existe $C>1$ tel que
$\ds{\frac{a}{C}\leq b \leq C \ a }$.)

\noindent {\bf  Remarque -} {\it Si $\eta \in V(o, \xi(t), t)$
l'angle en $\xi(t)$ entre les segments $[o, \xi(t)]$
et $[\xi(t), \eta)$ est sup\'erieur \`a $\pi/2$ ; par cons\'equent 
la diff\'erence $\cB_{\eta}(0, \xi(t))-t$ est born\'ee
uniform\'ement en $\eta \in V(o, \xi(t), t)$. En utilisant la 
$\delta$-conformit\'e de la famille de mesures $(\sigma_{x})_{x}$, le 
th\'eor\`eme pr\'ec\'edent peut donc s'\'enoncer de fa\c{c}on \'equivalente }
$$\sigma_{\xi(t)}(V(o, \xi, t))\asymp  
e^{(k_{t}-\delta)d(\xi(t), G.o)}.
$$

\noindent La d\'emonstration du th\'eor\`eme 5.1.  se fait en 
plusieurs \'etapes. Dans un premier temps nous consid\'erons 
le cas o\`u
$\xi_{t }\in \tC_{0}$ ; le 
th\'eor\`eme appara"t alors  comme une version am\'elior\'ee du th\'eor\`eme de 
l'ombre de Sullivan (lemme 2.2). D'apr\`es la remarque pr\'ec\'edente, il suffit de 
montrer que dans ce cas la quantit\'e $\sigma_{\xi(t)}(V(o, \xi, t))$ est 
comprise entre deux constantes strictement positives ind\'ependantes 
de $t$.

\noindent  Puisque $\xi(t)
$ appartient \`a $\tC_{0}$ il existe un point $g.o$ de $G.o$ \`a 
distance $\leq diam(C_{0})$ de $\xi(t)$. On obtient alors
$$
\sigma_{\xi(t)}(V(o, \xi, t)) \asymp  \sigma_{g.o}(V(o, \xi, t))$$
avec $\sigma_{g.o}(V(o, \xi, t))=
\sigma_{o}(V(g^{-1}.o, g^{-1}.\xi, t))$.  Ainsi 
$\sigma_{\xi(t)}(V(o, \xi, t))\leq \sigma _{o}(\SS^d)$. 
D'autre part, puisque $g^{-1}.\xi(t)$ est \`a distance $\leq  diam (C_{0})$
de $o$ on remarque qu'il existe 
$\epsilon >0$ tel que la boule euclidienne $B_{e}(g^{-1}.\xi, 
\epsilon)$ soit contenue dans 
$V(g^{-1}.o, g^{-1}.\xi, t)$. On conclut en notant que, pour $\epsilon > 0$
fix\'e,  on a 
$\ds{\inf_{\eta \in \Lambda_{G}}\sigma_{o}(B_{e}(\eta, 
\epsilon))>0}$.

\noindent Consid\'erons \`a pr\'esent le cas o\`u il existe $1\leq i \leq 
l$ tel que $\ds{\xi(t)\in
\tC_{i}}$. Soit $g$ l'\'el\'ement de $G$ tel que  
$\xi(t)\in g.\cH_{i}$ ; il existe alors  $0 \leq s \leq t$ tels que 
$[0, \xi(t)]\cap g. \cH_{i}=[\xi(s), \xi(t)]$. On a 
$\sigma_{o}(V(o, \xi, t)\asymp e^{-\delta s}
\sigma_{\xi(s)}(V(o, \xi, t))$ avec $d(\xi(s), \xi(t))= t-s$ et
$d(\xi(s), G.o) \leq diam(C_{0})$. Pour d\'emontrer le th\'eor\`eme 
il suffit de v\'erifier
que
$$
\sigma_{\xi(s)}(V(o, \xi, t))\asymp
\exp\bigl( -\delta(t-s)+(k_{i}-\delta)d(\xi(t), G.\xi(s))\bigr).
$$
On peut donc se ramener au cas o\`u $\xi(t) \in 
\cH_{i}$
 et $o \in \partial \cH_{i}$. \underline{Dans la suite de la d\'emonstration} nous 
 supposerons que le segment $[o, \xi(t)]$ est inclus dans ${\cal 
 H}_{i}$. 
 
\noindent Enon\c{c}ons d'abord le

\noindent {\bf Lemme 5.2.}
{\it Soit $\xi_i$  un point fixe parabolique born\'e de rang $k_{i}$ dans 
$\Lambda_{G}$. Alors }
$$
\sigma_{\xi_{i}(t)}(V(o, \xi_i, t)) \asymp e^{(k_{i}-\delta)t}
\quad \mbox{\rm et } \quad \sigma_{\xi_{i}(t)}(\SS^d- V(o, \xi_i, t)) 
\asymp e^{(k_{i}-\delta)t}.
$$

\noindent {\it D\'emonstration-}   Pla\c{c}ons nous 
dans le mod\`ele du demi-espace avec 
$\xi_i = +\infty$.
Dans ce mod\`ele, les points $o$  et $\xi_{i}(t)$ correspondent 
respectivement  aux couples $(0, 1)$ et 
$(0, e^t)$ et l'ensemble $\SS^d- V(o, \xi_i, t)$ n'est 
autre que le disque euclidien 
sur  $\C$  centr\'e en 
$0$ et de rayon $e^t$.
Notons $P$ le stabilisateur de $\xi_i$ dans $P$.
Puisque $G$ est g\'eo\-m\'e\-tri\-quement fini, on peut choisir
 un domaine fondamental bor\'elien $\cD$
 pour l'action de $ P$ sur $\C$
tel que $\cD$ contienne $0$ et tel que  
$  \cD\cap \Lambda_{G}$ soit relativement compact dans $\C$.

\noindent Le groupe $G$ \'etant divergent, on a
$\sigma_{o}(\xi_{i}) = 0$ d'o\`u
$$
 \sum_{p/ \vert o-p.o\vert 
\geq  e^t+\Delta }\sigma_{\xi_{i}(t)}(p \cD)
\leq \sigma_{\xi_{i}(t)}(V(o, \xi_i, t))\leq  \sum_{p/ \vert o-p.o\vert 
\geq  e^t-\Delta }\sigma_{\xi_{i}(t)}(p \cD)
$$
et 
$$
 \sum_{p/ \vert o-p.o\vert 
\leq  e^t-\Delta } \sigma_{\xi_{i}(t)}(p \cD)
\leq \sigma_{\xi_{i}(t)}(\SS^d-V(o, \xi_i, t))\leq  \sum_{p/ \vert o-p.o\vert 
\leq  e^t+\Delta }\sigma_{\xi_{i}(t)}(p \cD)
$$
o\`u   $\Delta$ d\'esigne le diam\`etre de $\cD$ dans $\C$.
Par convexit\'e des horisph\`eres, la fonction de Buseman
$\cB_{\eta}(\xi_{i}(t), p.\xi_{i}(t))$ diff\`ere de $d(\xi_{i}(t), 
p.\xi_{i}(t))$ 
d'une quantit\'e 
born\'ee uniform\'ement par 
rapport \`a
$p \in P$ et $\eta \in p\cD$ ; par cons\'equent 
$$
\sigma_{\xi_{i}(t)}(p \cD) \asymp e^{-\delta d(\xi_{i}(t), 
p.\xi_{i}(t))}\sigma_{p.\xi_{i}(t)}(p.\cD)
\asymp e^{-\delta t-\delta d(\xi_{i}(t), 
p.\xi_{i}(t))}\sigma_{o}(\cD)$$

\noindent Il nous reste \`a  estimer $d(\xi_{i}(t), 
p.\xi_{i}(t))$ ; on a 
$\ds{e^{-d(\xi_{i}(t), 
p.\xi_{i}(t))}= \frac{1-T}{1+T}}$ avec 
$\ds{T = \sqrt{\frac{\vert o-p.o\vert^{2}}
{\vert o-p.o\vert^{2}+4e^{2t}}}}$ 
d'o\`u $\ds{e^{-d(\xi_{i}(t), 
p.\xi_{i}(t))}\asymp \frac{1}{1+(\vert o-p.o\vert/e^t)^{2}}}$.

\noindent Si $P = < p_{\alpha}>$,  le point $\xi_{i}$ est de rang $1$; 
on  a
$$
\sigma_{\xi_{i}(t)}(V(x, \xi_i, t))\asymp e^{-\delta t} \sum_{n/ \vert n 
\alpha  \vert >e^t}\Bigl(\frac{1}{1+(\vert n\alpha \vert/e^{t})^{2}}\Bigr)^{\delta}
\asymp e^{-\delta t}\int_{e^t}^{+\infty}\Bigl(\frac{1}{1+(u/e^t)^{2}}\Bigr)^\delta du \asymp 
e^{(1-\delta) t}
$$
et  
$$
\sigma_{\xi_{i}(t)}(\SS^d- V(x, \xi_i, t))
\asymp e^{-\delta t}\int_{0}^{e^t}\Bigl(\frac{1}{1+(u/e^t)^{2}}\Bigr)^\delta du 
\asymp 
e^{(1-\delta) t}.
$$
Lorsque   $P = < p_{\alpha}, p_{\beta}>$, le point $\xi_{i}$ est de 
rang $2$; il vient
\begin{eqnarray*}
\sigma_{\xi_{i}(t)}(V(o, \xi_i, t)) &\asymp& e^{-\delta t} 
\sum_{n, m/ \vert n\alpha + m \beta  \vert >e^t}
\Bigl(\frac{1}{1+(\vert n\alpha + m \beta 
\vert/e^{t})^{2}}\Bigr)^{\delta}\\
&\asymp& e^{-\delta t}\int_{e^t}^{+\infty}\int_{0}^{2\pi}
\Bigl(\frac{1}{1+(r/e^t)^{2}}\Bigr)^\delta 
r\ dr \ d\theta  \asymp 
e^{(2-\delta) t}
\end{eqnarray*}
et  
$$
\sigma_{\xi_{i}(t)}(\SS^d- V(o, \xi_i, t))\asymp 
e^{-\delta t}\int_0^{e^t}\int_{0}^{2\pi} 
\Bigl(\frac{1}{1+(r/e^t)^{2}}\Bigr)^\delta 
r\ dr \ d\theta  \asymp 
e^{(2-\delta) t}. 
$$
\begin{flushright} \fdem
\end{flushright}

\noindent Le th\'eor\`eme 5.1.  est ainsi d\'emontr\'e lorsque $\xi$ est un 
point  parabolique et que $\xi(t)$ se projette sur la partie cuspidale 
correspondante. Cherchons maintenant \`a estimer $\sigma_{\xi(t)}(V(o, \xi, 
t))$ 
lorsque $\xi$ n'est plus un point parabolique mais 
que le point $\xi(t)$  
se projette n\'eanmoins sur la partie mince de $M$. 
 Dans la suite, on fixe une constante $C>0$.

 {\it Supposons dans un premier temps que le point 
$\xi_{i}$ appartient \`a $V(o, \xi, t+C)$}.

\noindent 
Il existe alors   $c>0$ (ne d\'ependant que de $C$)
tel  que
$$
V(o,\xi_i, t+c)\subset V(o, \xi, t)\subset V(o, \xi_i, t-c).
$$
(pour s'en convaincre, on peut se rappeler que
$V(o, \xi, t)$ est une boule euclidienne de $\SS^d$ de centre $\xi$ 
et de rayon $r_{t}\asymp e^{-t}$; on a alors $\vert \xi - 
\xi_{i}\vert <r_{t+C}$ et l'on peut choisir $\rho<1$
tel que
$B_{e}(\xi_{i}, \rho r_{t}) \subset B_{e}(\xi,  r_{t}) 
\subset B_{e}(\xi_{i}, r_{t}/\rho)$ puis  $c>0$ tel que 
$V(o,\xi_i, t+c) \subset B_{e}(\xi_{i}, \rho r_{t})$ et 
$B_{e}(\xi_{i}, r_{t}/\rho)\subset  V(o, \xi_i, t-c)$).

\noindent Remarquons que la distance entre $\xi(t)$ et $\xi_i(t)$ est major\'ee 
ind\'ependamment de $t$. En effet, un argument de convexit\'e
montre que l'angle entre 
les segments g\'eod\'esiques $[o, \xi(t)]$ et 
$[\xi(t), \xi_{i}(t+c)]$ est minor\'e par une constante 
strictement positive ; la diff\'erence entre 
$d(o, \xi(t))+d(\xi(t), \xi_i(t+c))$ et 
$d(o, \xi_{i}(t+c))$ est donc uniform\'ement born\'ee et il en est  de m\^eme
pour  $d(\xi_i(t+c), \xi(t))$. Par cons\'equent
$\ds{
\sigma_{\xi(t)}(V(o, \xi, t))\asymp \sigma_{\xi_i(t)}(V(o, \xi, t))
}.$

\noindent Par le lemme pr\'ec\'edent,  on sait que
$\ds{\sigma_{\xi_i(t)}(V(o, \xi_i, t+c))\asymp e^{(k-\delta)t}}$ ;
on conclut en remarquant que $t-d(\xi(t), G.o)$ est 
  uniform\'ement born\'e par rapport \`a $t$.

\noindent L'estimation de 
$\ds{
\sigma_{\xi(t)}(\SS^d-V(o, \xi, t))
}$ s'obtient de la m\^eme fa\c{c}on \`a partir de
la double inclusion
$
\SS^d-V(o,\xi_i, t-c)\subset \SS^d-V(o, \xi, t)\subset \SS^d-V(o, \xi_i, t+c).
$

 {\it  Supposons maintenant que }
$\xi_i \notin V(o, \xi, t-C)$.

\noindent Posons  $o' = \partial \cH_{i}\cap ]o, \xi[$ et  
$t'=d(\xi(t), o')$. Si  $\xi'$   d\'esigne le point antipodal de 
$\xi$, on voit que   
$$V(o, \xi, t)= \SS^d- V(o', \xi', t')
\quad \mbox{\rm et} \quad \xi_i \in V(o', \xi', t'+C).$$
Soit $p.o$ le point de l'orbite  $P.o$ le plus proche de $o'$ 
; la distance $d(p^{-1}.o', o)$ est major\'ee 
par $diam(C_{0})$
et 
$\xi_i$ appartient \`a
$V(p^{-1}.o', p^{-1}.\xi', t'+C).$ De plus, il existe $c'>0$ ne 
d\'ependant que de $C$ et de $diam(C_{0})$ tel que
$$
V(o, p^{-1}.\xi', t'-c')
\subset V(p^{-1}.o', p^{-1}.\xi', t'+C)
\subset V(o, p^{-1}.\xi', t'+c').
$$
Puisque le point $\xi_i$ appartient \`a
$V(o, p^{-1}.\xi', t'+c')$ on a  finalement 
\begin{eqnarray*}
\sigma_{\xi(t)}(V(o, \xi, t))&= &\sigma_{\xi(t)}(\SS^d- V(o', \xi', 
t'))\\
&=&\sigma_{p^{-1}\xi(t)}(\SS^d- V(p^{-1}.o', p^{-1}.\xi', 
t')\\
&\asymp&
\sigma_{x'_{t}}(\SS^d- V(o, p^{-1}.\xi', 
t')
\end{eqnarray*}
o\`u $x'_{t}$ d\'esigne le point de $[o, p^{-1}.\xi')$ situ\'e \`a 
distance $t'$ de $o$.
D'apr\`es les estimations ci-dessus, on a 
$$
\sigma_{x'_{t}}(\SS^d- V(o, p^{-1}.\xi', 
t') \asymp e^{(k_{i}-\delta)t'}
$$
et l'on conclut en   remarquant que $t'-d(\xi(t), G.o)$ est uniform\'ement 
born\'e.

{\it  Examinons enfin le cas o\`u
$\xi_i \notin V(o, \xi, t+C)$ mais $\xi_i \in V(o, \xi, t-C)$.}

\noindent Remarquons que 
$V(o, \xi, t+2C)\subset V(o, \xi, t)\subset V(o, \xi, t-2C)$.
En \'ecrivant alors $t-C=(t-2C)+C$ et $t+C=(t+2C)-C$, on obtient, gr\^ace 
aux
 deux \' etapes pr\'ec\'edentes
$$\sigma_{\xi(t-2C)}V(o, \xi, t-C) \asymp e^{(k_{i}-\delta)t}
\quad
\mbox{\rm et } \quad 
\sigma_{\xi(t+2C)}V(o, \xi, t+C) \asymp e^{(k_{i}-\delta)t}
$$
d'o\`u le r\'esultat.\fdem

\noindent Cette estimation de la mesure de Patterson des boules  de 
$\SS^d$ centr\'ees en des 
points de l'ensemble limite est essentielle pour \'etablir le r\'esultat 
suivant  du \`a D. Sullivan. 

\noindent {\bf Th\'eor\`eme 5.3.} 
{\it
Soit $G$ un groupe g\'eom\'etriquement fini et $k_{1}, \cdots 
k_{l}$ le rang des diff\'erents bouts cuspidaux de la vari\'et\'e 
quotient $M = \HH^{d+1}/G$. Notons $\sigma_{o}$ la mesure de Patterson de 
$G$ vue de l'origine, ${\cal H}_{\delta}$ la mesure de Hausdorff  et 
${\cal P}_{\delta}$ la mesure de packing de de $(\Lambda_{G}, \vert . 
\vert)$ 
associ\'ees \`a la jauge $r^\delta$.

Si pour tout $i\in \{1, \cdots, l\}$ on a $k_{i} \geq \delta$ alors 
$\sigma_{o}={\cal P}_{\delta}$. S'il existe $i \in \{1, \cdots, 
l\}$ tel que $k_{i}> \delta$  alors ${\cal 
H_{\delta}}$ est nulle sur tout  sous-ensemble de 
$\Lambda_{G}$ de $\sigma_o$-mesure positive.

Si pour tout $i\in \{1, \cdots, l\}$ on a $k_{i} \leq \delta$ alors $\sigma_{o}=
{\cal 
H_{\delta}}$. S'il existe $i \in \{1, \cdots, 
l\}$ tel que $k_{i} < \delta$  alors  ${\cal P}_{\delta}$ 
est infinie sur tout  sous-ensemble de 
$\Lambda_{G}$ de $\sigma_o$-mesure positive.}

\noindent {\it D\'emonstration-}
Supposons que pour tout $i\in \{1, \cdots, l\}$ on ait $k_{i} \geq \delta$. 
D'apr\`es les estimations pr\'ec\'edentes, il existe $c>0$ tel que 
pour tout point $\xi \in 
\Lambda_{G}$ on ait
$$\frac{\sigma_{o}(B_{e}(\xi, r))}{r^{\delta}}\geq c.$$
Par ailleurs, la mesure de Bowen-Margulis \'etant ergodique  pour le 
flot g\'eod\'esique sur $T^{1}M$ il vient
$$
\limsup_{t\to +\infty} {\bf 1}_{T^{1}C_{0}}(\phi_{t}v)=1 \quad 
\mu^{\sigma}(dv)-p.s.
$$
puisque $\mu^{\sigma}(T^{1}C_{0})>0$. Par cons\'equent 
$\ds{\limsup_{t \to +\infty} {\bf 1}_{C_{0}}(\xi(t))=1}\quad 
\sigma(d\xi)-p.s.$ et par le th\'eor\`eme pr\'ec\'edent il existe 
$C>0$ tel que pour $\sigma_o$-presque tout $\xi$ on ait
$$
\liminf_{r \to 0} \frac{\sigma_{o}(B_{e}(\xi, r))}{r^{\delta}} \leq C.
$$
Un argument de type Frostman montre qu'il existe alors $K>1$ tel que
$$\frac{1}{K} {\cal P}_{\delta} \leq \sigma_{o} \leq  K {\cal P}_{\delta}.$$
La mesure ${\cal P}_{\delta}$ \'etant une mesure conforme sur 
$\Lambda_{G}$ d'exposant $\delta$ elle coincide n\'ecessairement  
avec $\sigma_o$  par unicit\'e d'une telle mesure.

\noindent Supposons maintenant qu'il existe $k_{i}>\delta$ ; par 
ergodicit\'e de la mesure de Bowen-Margulis on a comme
pr\'ec\' edemment 
$$
\limsup_{t\to +\infty} d(\pi(\phi_{t}v), o) 
{\bf 1}_{T^{1}C_{i}}(\phi_{t}v)=+\infty \quad 
\mu^{\sigma}(dv)-p.s.
$$ Il vient  
$\ds{
\limsup_{r \to 0} \frac{\sigma_{o}(B_{e}(\xi, r))}{r^{\delta}} =+\infty
}\quad  \sigma_{o}(d\xi)-p.s.$
puisque $k_{i}>\delta$.    L'argument de Frostman permet de 
conclure.

\noindent Le cas 
o\`u tous les points paraboliques sont de rang $\leq \delta$ se traite 
de fa\c{c}on analogue.\fdem

\noindent {\bf R\'ef\'erences-}

\noindent Stratmann B.,  Velani S. L.
{\it The Patterson measure for geometrically finite groups with 
parabolic elements, new and old}
Proc. Lond. Math. Soc. 71 (1995) 197-220.

\noindent Sullivan D. 
{\it  Entropy, Hausdorff measures old and new, and limit sets 
of geometrically finite Kleinian groups}
Acta Math., (1984) p. 259-277.
\end{document}